\documentclass[a4paper]{article}

\usepackage{amsmath}
\usepackage{amssymb}

\numberwithin{equation}{section}

\newtheorem{thm}{Theorem}[section]
\newtheorem{prop}[thm]{Proposition}
\newtheorem{cor}[thm]{Corollary}
\newtheorem{dfn}[thm]{Definition}
\newtheorem{cjt}[thm]{Conjecture}
\newtheorem{lem}[thm]{Lemma}

\newcommand{\nn}{\nonumber}
\newcommand{\p}{\partial}
\newcommand{\ve}{\epsilon}
\newcommand{\dl}{\delta}
\newcommand{\LM}{{\mathcal L}(M)}
\newcommand{\Ker}{\mathrm{Ker}}
\newcommand{\Img}{\mathrm{Im}}

\newcommand{\pal}{\partial}

\newcommand{\al}{\alpha}

\newenvironment{prf}{\noindent {\it Proof} \ }{\hfill $\Box$}
\newenvironment{prfM}{\noindent {\it Proof of Theorem \ref{mthm}} \ }{\hfill $\Box$}
\begin{document}

\title{Deformations of Semisimple Bihamiltonian Structures of
Hydrodynamic Type}
\author{{Si-Qi Liu \ \ Youjin Zhang}\\
{\small Department of Mathematical Sciences, Tsinghua
University}\\
{\small Beijing 100084, P.R.China}\\
{\small  lsq99@mails.tsinghua.edu.cn,\
yzhang@math.tsinghua.edu.cn}}
\date{}
\maketitle
\begin{abstract}
We classify in this paper infinitesimal quasitrivial deformations of  semisimple
bihamiltonian structures of hydrodynamic type.
\end{abstract}

\section{Introduction}
A bihamiltonian structure of hydrodynamic type defined on the formal loop space of a manifold $M$
consists of two compatible Poisson brackets of the form
\begin{equation}\label{z3}
\{u^i(x),u^j(y)\}=g^{ij}(u(x))\delta'(x-y)+\Gamma^{ij}_{k}(u(x))
u^k_x \delta(x-y),\quad i,j=1,\dots, n.
\end{equation}
Here $n=\dim M$, and we assume that $\det(g^{ij}(u))\ne 0$. Such type of Poisson brackets were
introduced and classified by Dubrovin and Novikov during the 80's
of the last century \cite{dn83, dn84,dn89}, they were used to describe the hamiltonian structures of systems of
hydrodynamic type. According to the theory of Dubrovin and Novikov, the inverse of $(g^{ij})$
must be a flat metric of the manifold $M$, and the coefficients $\Gamma^{ij}_k$ be given by the contravariant
components of the Levi-Civita connection of this flat metric. Two such Poisson brackets corresponding
to two flat metrics $(g^{ij}_1)^{-1}, (g^{ij}_2)^{-1}$ are compatible if these two metrics form a flat pencil \cite{Du1}.
The most well known examples of bihamiltonian structures of hydrodynamic type
are possessed by the Whitham equations (in particular, the dispersionless limit) of
integrable evolutionary PDEs of KdV type \cite{dn83, dn84,dn89,maltsev}.

In the present paper we study the problem of classification of deformations of a given
bihamiltonian structure of hydrodynamic type, these deformations depend on a parameter $\ve$ which
is called the dispersion parameter. The deformed bihamiltonian structure has the form
\begin{eqnarray}
&&\{u^i(x),u^j(y)\}_{a}=g^{ij}_a(u(x))\delta'(x-y)+\Gamma^{ij}_{k;a}(u(x)) u^k_x \delta(x-y)
\nn\\&&+\sum_{m\ge 1}\sum_{l=0}^{m+1} \ve^m A^{ij}_{m,l;a}(u;u_x,\dots,u^{(m+1-l)}) \delta^{(l)}(x-y),
\quad  a=1,2.\label{z1}
\end{eqnarray}
Here $A^{ij}_{m,l;a}$ are {\em differential polynomials}, i.e. they depend polynomially
on the $x$-derivatives of $u^1,\dots,u^n$, and the
coefficients of these polynomials are smooth functions of $u^1,\dots, u^n$.
We also require that $A^{ij}_{m,l;a}$ are homogeneous polynomials in the sense that if we assign degree
$m$ to $u^{i,m}=\p^m_x u^i$, then $\deg A^{ij}_{m,l;a}=m+1-l$.
The class of bihamiltonian structures of the form (\ref{z1}) that satisfy some additional
conditions is classified in \cite{DZ1}.
These additional conditions include the so called tau-symmetry property and the linearization of the
Virasoro symmetries of the corresponding hierarchy of bihamiltonian evolutionary PDEs, they ensure the existence
of tau functions for solutions of the hierarchy and the possibility of
representing the Virasoro symmetries of the hierarchy by the action of an infinite
number of linear differential operators on the tau functions. The moduli
space of this class of bihamiltonian structures coincides with the space of semisimple Frobenius manifolds\cite{DZ1}.
Here we will study the class of deformed bihamiltonian structures of the form (\ref{z1})
without the restriction of these additional properties.

The bihamiltonian structures of hydrodynamic type under our considerations are assumed to be
semisimple,
i.e., the eigenvalues of the matrix $({g_1}^{ij})^{-1} g_2^{ij}$ are pairwise distinct, here
$({g_1}^{ij})^{-1}, (g_2^{ij})^{-1}$ are the flat metrics corresponding to the given bihamiltonian structure.
The simplest example of semisimple bihamiltonian structures of hydrodynamic type has the form
\begin{eqnarray}
&&\{u(x),u(y)\}_1=\delta'(x-y),\nn\\
&&\{u(x),u(y)\}_2=u(x)\delta'(x-y)+\frac12\,u(x)'\delta(x-y),\label{z2}
\end{eqnarray}
it is the dispersionless limit of the bihamiltonian structure of the KdV hierarchy \cite{gardner,magri,ZF}.
In \cite{Lo} Lorenzoni studied its deformations at the approximation up to $\ve^4$. He showed that the equivalence
classes of all such deformations are parameterized by a smooth function $s(u)$, the bihamiltonian structure
of the KdV hierarchy corresponds to the special deformation with a nonzero constant $s(u)$.
Here the equivalence relation between deformations of a bihamiltonian structures of hydrodynamic type
is defined in \cite{DZ1}, two deformations
of the form (\ref{z1}) are defined to be equivalent if they are related
by a Miura-type transformation
\begin{equation}\label{miura}
u^i\mapsto u^i+\sum_{k\ge 1}\ve^k F^i_k(u;u_x,\dots,u^{(k)}),\quad i=1,\dots,n
\end{equation}
where $F^i_k$ are {\em differential polynomials} of degree $k$, note that they are not required to depend polynomially
on $u^1,\dots, u^n$. In
particular, a deformation (\ref{z1}) is called to be trivial if it
is equivalent to the undeformed bihamiltonian structure. For the
above example, when the function $s(u)$ does not vanish, the
corresponding deformation of the bihamiltonian structure
(\ref{z2}) is nontrivial. Nevertheless, Lorenzoni proved that at
the approximation up to $\ve^4$ all such deformations are {\em
quasitrivial}. The notion of quasitriviality was also
introduced in \cite{DZ1}, a bihamiltonian structure of the form
(\ref{z1}) is called quasitrivial if it can be obtained from its
dispersionless limit by a transformation of the form
\begin{equation}\label{qmiura}
u^i\mapsto u^i+\sum_{k\ge 1} \ve^k G^i_k(u;u_x,\dots,u^{(m_k)}),\quad i=1,\dots,n.
\end{equation}
Here $G^i_k$ are smooth functions of their arguments, in particular,
they are not necessary polynomials of the $x$-derivatives of $u^1,\dots, u^n$.
In \cite{DZ1} it was proved that all semisimple bihamiltonian structures of the form (\ref{z1}) that
satisfy the tau-symmetry property are quasitrivial. The method given in there can in fact
be employed to prove the quasitriviality of all deformations of (\ref{z2}). These results suggest
that quasitriviality could hold true for any deformation (\ref{z1}) of a semisimple bihamiltonian structure of
hydrodynamic type.

In this paper we will restrict ourselves to study properties of quasitrivial deformations
and leave the discussion on the validity of quasitriviality for any deformation of a semisimple
bihamiltonian structure of hydrodynamic type to a subsequent publication.
The main result of the paper is contained in the following two theorems:
\begin{thm}\label{mthm-z1}
Any two quasitrivial deformations of a semisimple bihamiltonian structure of
hydrodynamic type are equivalent if and only if they are equivalent at the approximation
up to $\ve^2$.
\end{thm}
The semisimplicity of a bihamiltonian structure of hydrodynamic type implies the existence of a coordinate
system under which the corresponding two flat metrics are diagonal \cite{FEV},
we call such coordinates the canonical
coordinates of the semisimple  bihamiltonian structure.
\begin{thm}\label{mthm-z2}
At the approximation up to $\ve^2$, the space of the equivalence classes of all quasitrivial deformations
of a semisimple bihamiltonian structure of hydrodynamic type is parameterized by $n$ smooth
functions $c_1(u^1), \dots, c_n(u^n)$ of its canonical coordinates.
\end{thm}
We will prove the above theorems by classifying the infinitesimal quasitrivial deformations of
a given semisimple bihamiltonian structure of hydrodynamic type, it amounts to the calculation of
certain modification of the second bihamiltonian cohomology.
As a direct consequence of the calculation that will be performed
in section \ref{sec-4}, we have
\begin{cor} The equivalence classes of infinitesimal quasitrivial deformations of a semisimple bihamiltonian
structure of hydrodynamic type are parameterized by $n$ arbitrary
functions of one variable.
\end{cor}
The notion of bihamiltonian cohomology
was introduced in \cite{DZ1}, it provides an efficient tool to
study deformations of bihamiltonian structures.
We will first recall the notions of Poisson cohomology and bihamiltonian cohomology in section \ref{sec-2} and
section \ref{sec-3} respectively, and then give the proof of the main results in section \ref{sec-4},
some examples will be given in section \ref{sec-5}.

\section{Local Poisson structures and Poisson cohomologies}\label{sec-2}

We recall in this section the definition of local Poisson structures and Poisson cohomologies
that was presented in \cite{DZ1} on the formal loop space
$
{\mathcal L}(M)=\{S^1\to M\}
$
of a manifold $M$ of dimension $n$, we will closely follow the notations of \cite{DZ1}.
Choose a
chart $U$ on $M$ with local coordinates $u^1,\dots, u^n$, we denote
by ${\mathcal A}={\mathcal A}(U)$ the ring of differential polynomials of the form
$$
f(x,u,u_x,\dots)=\sum_{i_1,s_1,\dots,i_m,s_m} f_{i_1,s_1;\dots;i_m,s_m}(x;u)
u^{i_1,s_1}\dots u^{i_m,s_m}.
$$
Here $u=(u^1,\dots, u^n),\ u^{(s)}=(u^{1,s},\dots,u^{n,s})$ with $u^{i,s}=\frac{d^s u^i(x)}{d x^s}$,
and the coefficients of these differential polynomials
are smooth functions on $S^1\times M$. Denote
$$
{\mathcal A}_{0}={\mathcal A}/{\mathbb R},\quad {\mathcal A}_{1}={\mathcal A}_{0} dx,\quad
\Lambda={\mathcal A}_{1}/d {\mathcal A}_{0}
$$
where the operator $d: {\mathcal A}_{0}\to {\mathcal A}_{1}$ is defined by
$$
f\mapsto d f=\left(\frac{\pal f}{\pal x}+\sum \frac{\pal f}{\pal u^{i,s}} u^{i,s+1}\right)
dx.
$$
Elements of $\Lambda$ are called local functionals on ${\mathcal L}(M)$, they will be expressed
as integrals over $S^1$ of a representative differential polynomial
\begin{equation}
\int f(x;u(x),u_x(x),\dots,u^{(N)}(x)) dx.
\end{equation}
Later in Section \ref{sec-4} we will also use functionals of the above form with densities
$f$ being smooth functions of their arguments instead of being differential polynomials.

A local $k$-vector on the formal loop space is defined to be a formal infinite sum of the following form
\begin{equation}
\al=\sum \frac1{k!}\pal_{x_1}^{s_1}\dots\pal_{x_k}^{s_k} A^{i_1\dots i_k}
\frac{\pal}{\pal u^{i_1,s_1}(x_1)}\wedge\dots \wedge\frac{\pal}{\pal u^{i_k,s_k}(x_k)}
\end{equation}
with the coefficients $A$'s having the form
\begin{equation}
A^{i_1\dots i_k}=\sum_{p_2,\dots,p_{k}\ge 0} B^{i_1\dots i_k}_{p_2\dots p_k}(u(x_1);
u_x(x_1),\dots) \delta^{(p_2)}(x_1-x_2)\dots \delta^{(p_k)}(x_1-x_k).
\end{equation}
Here $B^{i_1\dots i_k}_{p_2\dots p_k}(u(x_1);
u_x(x_1),\dots)\in{\mathcal A}$, and
\begin{equation}
A^{i_1\dots i_k}=A^{i_1\dots i_k}(x_1,\dots,x_k;u(x_1),\dots,u(x_k),\dots)
\end{equation}
are antisymmetric with respect to the simultaneous permutations
$
i_p,x_p\leftrightarrow i_q, x_q.
$
These coefficients $A^{i_1\dots i_k}$ are called the components of the local $k$-vector $\al$.
The space of all such local $k$-vectors is denoted by $\Lambda_{loc}^k$.
In particular, a local vector field
on the formal loop space has the form
\begin{equation}\label{def-vec}
\xi=\sum_{i=1}^n\sum_{s\ge 0} \pal_x^s \xi^{i}(u(x);u_x(x),\dots)\frac{\pal}{\pal u^{i,s}(x)}
\end{equation}
which is also called a translation (along $x$) invariant evolutionary vector field. A local bivector takes
the form
\begin{equation}\label{bi-vec}
\omega=\frac12\sum \pal_{x}^{s} \pal_y^t \omega^{ij}\frac{\pal}{\pal u^{i,s}(x)}\wedge\frac{\pal}{\pal u^{j,t}(y)}
\end{equation}
with
\begin{equation}
\omega^{ij}=A^{ij}(x-y;u(x),u_x(x),\dots)=\sum_{k\ge 0} A^{ij}_k(u(x);u_x(x),\dots) \delta^{(k)}(x-y).
\end{equation}
It is assumed that the space $\Lambda_{loc}^0$ is the subspace of $\Lambda$ that consists of
local functionals of the form
\begin{equation}\label{def-func}
{\bar f}=\int f(u(x);u_x(x),\dots) dx,\quad f(u(x);u_x(x),\dots)\in {\mathcal A}_{0}.
\end{equation}

On the space of local multi-vectors
\begin{equation}
\Lambda^*_{loc}=\Lambda^0_{loc}\oplus\Lambda_{loc}^1\oplus\Lambda_{loc}^2\oplus\dots
\end{equation}
there is defined a bilinear operation of Schouten-Nijenhuis bracket
\begin{equation}
[\ ,\, ]: \ \Lambda^k_{loc}\times \Lambda^{l}_{loc}\to\Lambda^{k+l-1}_{loc},\quad k,l\ge 0
\end{equation}
By definition, the Schouten-Nijenhuis bracket of any two elements of $\Lambda^0_{loc}$ is equal to zero,
and the Schouten-Nijenhuis bracket of a local vector field $\xi$ of the form (\ref{def-vec})
with a local functional $\bar f$ of the form (\ref{def-func}) is defined by
\begin{equation}
[\xi,\bar f]=\int \sum\left(\pal_x^s \xi^i\right) \frac{\pal f(u(x);u_x(x),\dots)}{\pal u^{i,s}} dx
=\int \sum_{i=1}^n \xi^i \frac{\delta \bar f}{\delta u^i(x)} dx
\end{equation}
where
\begin{equation}
\frac{\delta\bar f}{\delta u^i(x)}=\sum_{s\ge 0} (-1)^s \pal_x^s\left(\frac{\pal f}{\pal u^{i,s}}\right).
\end{equation}
The Schouten-Nijenhuis bracket of two local vector fields is given by their usual commutator
\begin{eqnarray}
[\xi,\eta]&=&\sum \left(\xi^{j,t} \frac{\pal\eta^{i,s}}{\pal u^{j,t}}-\eta^{j,t}
\frac{\pal\xi^{i,s}}{\pal u^{j,t}}\right)\frac{\pal}{\pal u^{i,s}}\nn\\
&=&\sum \pal_x^s\left(\xi^{j,t} \frac{\pal\eta^{i}}{\pal u^{j,t}}-\eta^{j,t}
\frac{\pal\xi^{i}}{\pal u^{j,t}}\right)\frac{\pal}{\pal u^{i,s}},
\end{eqnarray}
and components of the Schouten-Nijenhuis bracket of a bivector $\omega$ of the form (\ref{bi-vec})
with a functional $I$ and with a local vector filed $\xi$ of the form (\ref{def-vec})
are given respectively by
\begin{eqnarray}
&& [\omega,I]^i=\sum_{j,k} A^{ij}_k \pal_x^k\frac{\delta I}{\delta u^j(x)},\\
&&[\omega, \xi]^{ij}=\sum_{k,t}\left(\pal_x^t \xi^k(u(x); \dots)
\frac{\pal A^{ij}}{\pal u^{k,t}(x)}
-\frac{\pal \xi^i(u(x); \dots)}{\pal u^{k,t}(x)} \pal_x^t A^{kj}\right.\nn\\
&&\quad \quad \qquad\left. -\frac{\pal
\xi^j(u(y); \dots)}{\pal u^{k,t}(y)} \pal_y^t A^{ik}\right).
\label{lie-der}
\end{eqnarray}

The Schouten-Nijenhuis bracket satisfies the following graded Jacobi identity and the antisymmetry property:
\begin{eqnarray}
&&(-1)^{k m}[[a,b],c]+(-1)^{kl}[[b,c],a]+(-1)^{lm} [[c,a],b]=0,\label{jacobi}\\
&&[a,b]=(-1)^{kl} [b,a],\quad a\in\Lambda_{loc}^k,\ b\in\Lambda_{loc}^l,\ c\in\Lambda_{loc}^m.
\end{eqnarray}

\begin{dfn}[\cite{DZ1}] A local bivector $\omega\in\Lambda^2_{loc}$ of the form
(\ref{bi-vec}) is called a local Poisson structure on the formal loop space ${\mathcal L}(M)$
if $[\omega,\omega]=0$.
\end{dfn}
A local Poisson structure given by a bivector of the form (\ref{bi-vec}) can also be represented
as an antisymmetric bilinear map from $\Lambda^2$ to $\Lambda$ as follows:
\begin{equation}
\{{\bar f_1},{\bar f_1}\}=\int \sum_{k\ge 0} \frac{\delta {\bar f_1}}{\delta u^i(x)}
A^{ij}_k(u(x);u_x(x),\dots)\pal_x^k \frac{\delta {\bar f_2}}{\delta u^j(x)} dx.
\end{equation}
For a particular choice of the local functionals ${\bar f_1}=\int u^i(z) \delta(z-x) dz,\
{\bar f_2}=\int u^j(z) \delta(z-y) dz$ we get the usual representation of a Poisson structure
\begin{equation}
\{u^i(x),u^j(y)\}=
\sum_{k\ge 0} A^{ij}_k(u(x);u_x(x),\dots) \delta^{(k)}(x-y).
\end{equation}

There is a natural gradation on the space of local multi-vectors which is defined by
\begin{equation}
\deg u^{i,s}=s,\ \deg\frac{\pal}{\pal u^{i,s}}=-s,\ \deg dx=-1,\ \deg\delta^{(s)}(x-y)=s+1.
\end{equation}
To separate monomials of different degree in a local multi-vector, we introduce a formal
indeterminate $\ve$ and assign to it the degree $-1$. Denote
\begin{eqnarray}
&&\Omega^k_m=\{a\in \left. \Lambda^k_{loc} \right|\deg a=m\},\nn\\
&&\Omega^k=\{a\in \left. \Lambda^k_{loc}\otimes{\mathbb{C}}[[\ve],\ve^{-1}] \right|\deg a=k\}.
\end{eqnarray}
For example, an element of $\Omega^0$ has the form
\begin{equation}
{\bar f}=\int \left({\ve}^{-1} f_0(u(x))+\sum_{k=1}^n f_{1,k}(u(x)) u^k_x+\dots\right) dx.
\end{equation}
The components of a vector field $\xi\in \Omega^1$ has the form
\begin{equation}
\xi^i=\ve^{-1} a^i(u)+\sum_{k=1}^n b^i_{k}(u) u^k_x+\ve
\left(\sum_{k=1}^n c^i_{k}(u) u^k_{xx}+ \sum_{k,l=1}^n e^i_{kl}(u) u^k_x u^l_x\right)+\dots.
\end{equation}
A Poisson structure $\omega\in \Omega^2_2$ is of hydrodynamic type and has the representation
of the form (\ref{z3}), any Poisson structure of the form $\omega+P(\ve)\in \Omega^2$ with
$P(\ve)=\sum_{k\ge 1}\ve^k P_k,\ P_k\in \Omega^2_{k+2}$ is called a deformation of $\omega$.

The space
\begin{equation}
\Omega=\Omega^0\oplus\Omega^1\oplus\Omega^2\oplus\dots
\end{equation}
is closed with respect to the Schouten-Nijenhuis bracket $\ve [\,,\,]$,
so a Poisson structure $\omega\in \Omega^2$ defines a differential
\begin{equation}
\ve d: \Omega^k\to \Omega^{k+1},\quad \ve\, d a=\ve [\omega,a],\quad  a\in \Omega^k.\label{df-p}
\end{equation}
The cohomology of the complex $(\Omega, \ve d)$ is called the Poisson cohomology of the Poisson
structure $\omega$, and is denoted by $H^*({\mathcal L}(M), \omega)$ \cite{DZ1}. It is a natural generalization of the
the notion of Poisson cohomology for finite dimensional Poisson structures \cite{lichn}.

\section{Bihamiltonian structures and bihamiltonian cohomologies}\label{sec-3}

Assume that we are given two Poisson structures $\omega_1, \omega_2$ of hydrodynamic type
with components of the form
\begin{equation}
\omega^{ij}_{a}=g^{ij}_{a}(u)\dl'(x-y)+\Gamma^{ij}_{k,a}(u)
u^k_x\dl(x-y),\quad \det(g^{ij}_a)\ne 0,\quad a=1,2.\label{ssbh}
\end{equation}
If the linear combination $\omega_{\lambda}=\omega_2-\lambda \omega_1$ is also a Poisson structure
for an arbitrary parameter $\lambda\in \mathbb{R}$, then the pair $(\omega_1,\omega_2)$ is called a
bihamiltonian structure of hydrodynamic type. These two Poisson structures
define two complexes $(\Omega,\ve d_a),\ a=1,2$. It is proved in \cite{magri2,get}
that the Poisson cohomologies $H^*(\LM,\omega_a), a=1,2$
are trivial (also see \cite{DZ1} for a different proof of triviality
for the first and the second Poisson cohomologies). Thus any deformation $\omega_a+P(\ve) \in \Omega^2$ of a
single Poisson structure $\omega_a$ can be obtained from $\omega_a$
by performing a Miura type transformation of the form (\ref{miura}).
Instead of the deformations of a single Hamiltonian structure, we are interested in deformations of
the bihamiltonian structure $(\omega_1, \omega_2)$. Due to the triviality of the Poisson cohomology
$H^*(\LM,\omega_1)$, we can always assume that our deformations keep the first
Poisson structure $\omega_1$ unchanged.

\begin{dfn} The pair of bivectors
\begin{equation}\label{def-bs}
(\omega_1,\ \omega_2+\sum_{m\ge 1}\ve^m P_m),\ P_m \in \Omega^2_{m+2}.
\end{equation}
is called a deformation of the bihamiltonian structure $(\omega_1,\omega_2)$
if the equality
\begin{equation} \label{def-sb}
[\omega_2+\sum_{m\ge 1}\ve^m P_m-\lambda\omega_1,\
\omega_2+\sum_{m\ge 1}\ve^m P_m-\lambda\omega_1]=0
\end{equation}
holds true for an arbitrary parameter $\lambda$. It is called an $N$-th order
deformation of the bihamiltonian structure $(\omega_1,\omega_2)$ if the
equlity (\ref{def-sb}) holds true for an arbitrary parameter $\lambda$ at the
approximation up to $\ve^N$.
\end{dfn}

\begin{dfn}
We say that two deformations (of order $N$) of the bihamiltonian structure
$(\omega_1,\omega_2)$ are equivalent or quasi-equivalent if they are related
(resp. at the approximation up to $\ve^N$) by a Miura type transformation (\ref{miura}) or
by a quasi-Miura type transformation (\ref{qmiura}). A deformation (of order $N$) of the
bihamiltonian structure $(\omega_1,\omega_2)$ is called trivial or quasitrivial if it is equivalent
or quasi-equivalent to $(\omega_1,\omega_2)$ (resp. at the approximation up to $\ve^N$).
\end{dfn}
Due to the above definition, for a $N$-th order deformation (\ref{def-bs})
the bivectors $P_m$ must satisfy the conditions
\begin{eqnarray}
&&d_1 P_m=0,\ 1\le m\le N,\label{SZ-3a}\\
&&d_2 P_1=0,\ 2\,d_2 P_{m}+\sum_{k=1}^{m-1} [P_k,P_{m-k}]=0,
\ 2\le m\le N.\label{SZ-3b}
\end{eqnarray}
Here the differentials $d_1, d_2$ are defined by the Poisson structures $\omega_1$ and $\omega_2$
respectively
as in (\ref{df-p}), they act on the subspaces $\Omega^k_m$ as
\begin{equation}
d_a: \Omega^k_m\to \Omega^{k+1}_{m+2},\ k\ge 0,\ m\ge k-1,\ a=1,2.
\end{equation}
The notion of bihamiltonian cohomologies
$
H^k=\oplus_{m\ge k-1} H^k_m,\, k\ge 0
$
for $(\omega_1, \omega_2)$ is introduced in \cite{DZ1}, they are defined by
\begin{eqnarray}
&&H^k_m(\LM;\omega_1,\omega_2)=\Ker(d_1 d_2|_{\Omega^{k-1}_{m}})/
\Img(d_1|_{\Omega^{k-2}_{m-2}})\oplus\Img(d_2|_{\Omega^{k-2}_{m-2}}),\quad k\ge 2,\nn\\
&&H^1_m(\LM;\omega_1,\omega_2)=\Ker(d_1 d_2|_{\Omega^{0}_{m}})\nn\\
&&H^0_m(\LM;\omega_1,\omega_2)=\Ker(d_1|_{\Omega^{0}_{m}})\cap \Ker(d_2|_{\Omega^{0}_{m}})
\end{eqnarray}
It was proved in \cite{DZ1} that the zero-th cohomology coincides with the space of common Casimirs of the
Poisson structures $\omega_1, \omega_2$, the first cohomology corresponds to the space of
bihamiltonian vector fields, and the second cohomology corresponds to the space of infinitesimal
deformations of the bihamiltonian structure modulo the trivial deformations caused by Miura transformations.
Below we list some other simple propositions on the second and third cohomologies.

\begin{prop}\label{prop-z1}
1). The bihamiltonian cohomologies $H^2_i(\LM;\omega_1,\omega_2)$
vanish for $K+1\le i\le N$ iff any class of deformations of the bihamiltonian
structure $(\omega_1,\omega_2)$ of order $s\le N$ is uniquely determined by the
corresponding class of deformations of order $K$; 2).
The bihamiltonian cohomologies $H^2_{2k+1}(\LM;\omega_1,\omega_2)$
vanish for $1\le 2k+1\le N$ iff any deformation of the bihamiltonian structure
$(\omega_1,\omega_2)$ is equivalent to a deformation of the form
(\ref{def-bs}) with $P_{2l+1}=0,\ 2 l+1\le N$.
\end{prop}
\begin{prf}
Let us first assume that $H^2_i(\LM;\omega_1,\omega_2)$ vanishes for
$K+1\le i\le N$. We need to prove that any two deformations of order $s\le N$ of the form
\begin{equation}\label{SZ-6}
(\omega_1,\,\omega_2+\sum_{m=1}^K\ve^m P_m+\sum_{m=K+1}^s\ve^m P^{(l)}_m)+
{\mathcal O}(\ve^{s+1}),\quad l=1,2
\end{equation}
are equivalent. By using the identities in (\ref{SZ-3a}),(\ref{SZ-3b}) we can find
$X,Y\in \Omega^1_{K+1}$ such that
$$P_{K+1}^{(1)}=d_1 X,\quad P_{K+1}^{(2)}=d_1 Y.$$
From (\ref{SZ-3b}) it follows that
$$d_2 d_1 (X-Y)=0.$$
So our assumption implies the existence of $I, J\in \Omega^0_{K-1}$ such that
$$X=Y+d_1 I+d_2 J.$$
Thus after the Miura type transformation
$$u^i\mapsto u^i-\ve^{K+1} d_1 J$$
the first deformation
$$(\omega_1,\,\omega_2+\sum_{m=1}^K\ve^m P_m+\sum_{m=K+1}^s\ve^m P^{(1)}_m)
+{\mathcal O}(\ve^{s+1})$$
is transformed to
$$(\omega_1,\,\omega_2+\sum_{m=1}^K\ve^m P_m+\ve^{K+1} P_{K+1}^{(2)}
+\sum_{m=K+2}^s\ve^m {\tilde P}^{(1)}_m)+{\mathcal O}(\ve^{s+1})$$
By repeating the same procedure, we prove the equivalence of the two
deformations of (\ref{SZ-6}).

Now we assume that any class of deformations of the bihamiltonian structure
$(\omega_1,\omega_2)$ of order $s\le N$ is uniquely determined by the
corresponding class of deformations of order $K$. For any
$$X\in \Ker(d_1d_2|_{\Omega^1_s}),\ K+1\le s\le N$$
we have a $s$-th order deformation of the form
\begin{equation}\label{SZ-7}
(\omega_1,\,\omega_2+\ve^s d_1 X).
\end{equation}
It follows from our assumption that there exists a Miura type transformation
$$
u^i\mapsto u^i+\sum_{j=1}^s \ve^j A_j^i,\quad A_j\in \Omega^1_j
$$
that transforms the bihamiltonian structure $(\omega_1, \omega_2)$ to
(\ref{SZ-7}), i.e.,
\begin{eqnarray}
&&\omega_1=e^{-\ve^s {ad}_{{\tilde A}_s}}\dots e^{-\ve {ad}_{{\tilde A}_1}} \omega_1+{\mathcal O}(\ve^{s+1}),\nn\\
&&
\omega_2+\ve^s d_1 X=e^{-\ve^s {ad}_{{\tilde A}_s}}\dots e^{-\ve {ad}_{{\tilde A}_1}}
\omega_2+{\mathcal O}(\ve^{s+1}).\label{z5}
\end{eqnarray}
Here we represent, modulo $\ve^{s+1}$, the Miura transformation
as the composition of the one parameter transformation
groups $u\mapsto e^{\ve^k {\tilde A}_k} u,\ k=1,\dots, s$ corresponding to the vector fields
$$
{\tilde A^i}_1=A^i_1,\
{\tilde A^i_2}=A^i_2-\frac12 \sum_{j=1}^n \sum_{t\ge 0} \frac{\p A_1^i}{\p u^{j,t}}\, \p_x^t A_1^j,\dots
$$
From the identities in (\ref{z5}) we obtain
$$d_1 {\tilde A}_s=0,\quad d_2 {\tilde A}_s=d_1 X.$$
The first equality yields the existence of $I\in\Omega^0_{s-2}$ such that
${\tilde A}_s=d_1 I$, and from the second equality it follows that
$X\in \Img(d_1|_{\Omega^0_{s-2}})\oplus\Img(d_2|_{\Omega^0_{s-2}})$.
Thus we proved the first part of the proposition. The second part
can be proved in a similar way. The proposition is proved.\end{prf}

\begin{prop}
If the bihamiltonian cohomology $H^3_{N+3}(\LM;\omega_1,\omega_2)$ vanishes
then any $N$-th order deformation of the bihamiltonian structure
$(\omega_1, \omega_2)$ can be extended to a $N+1$-th order deformation.
\end{prop}
\begin{prf}
Any $N$-th order deformation can be represented as
$$(\omega_1,\,\omega_2+\sum_{i=1}^N \ve^i d_1 X_i)+{\mathcal O}(\ve^{N+1}),
\quad X_i\in \Omega^1_i.$$
In order to extend it to a deformation of order $N+1$ we need to find a local vector
field
$X_{N+1}\in \Omega^1_{N+1}$
such that
\begin{equation}\label{SZ-8}
d_1 d_2 X_{N+1}=\frac12\,\sum_{i=1}^N [d_1 X_i, d_1 X_{N+1-i}].
\end{equation}
Denote by $Q$ the r.h.s. of the above equation. Then by using the graded Jacobi
identity (\ref{jacobi}) of the Schouten-Nijenhuis bracket and the equalities
$$d_1 d_2 X_m=\frac12 \sum_{i=1}^{m-1} [d_1 X_i, d_1 X_{m-i}],\quad m=1,\dots,N$$
we obtain
$$d_1 Q=d_2 Q=0.$$
So there exists $R\in \Omega^2_{N+3}$ such that $Q=d_1 R$. Now it
follows from the equality $d_1 d_2 R=0$ and our assumption of the proposition that
$$R=d_1 A+d_2 B,\quad A, B\in\Omega^1_{N+1}.$$
So the equation (\ref{SZ-8}) now takes the form
$$d_1 d_2 X_{N+1}=d_1 (d_1 A+ d_2 B)$$
and it has a solution $X_{N+1}=B$. The proposition is proved.\end{prf}

Due to the the above propositions, the problem of classification of deformations of
the hydrodynamic bihamiltonian structures is reduced to the computation of
bihamiltonian cohomology. We can also consider certain modification of the bihamiltonian
cohomology in order to deal with quasitrivial deformations of the hydrodynamic bihamiltonian structures,
we will do this in the next section.

\section{Computation of a modified bihamiltonian cohomology and the proof of the main theorems}\label{sec-4}
We consider in this section the problem of classification of infinitesimal quasitrivial deformations of a
semisimple bihamiltonian structure $(\omega_1, \omega_2)$
with components of the form
(\ref{ssbh}). Let us choose the coordinates $u^1,\dots,u^n$, called the canonical coordinates of the semisimple
bihamiltonian structure, such that both metrics
$g^{ij}_1$ and $g^{ij}_2$ are diagonal under these coordinates, and the identities $g^{ii}_2=u^i\,g^{ii}_1$
hold true \cite{FEV}.
In terms of these coordinates the
bihamiltonian structure can be expressed as
\begin{eqnarray}
&&\omega^{ij}_1=f^i\,\delta^{ij}\delta'(x-y)+\frac12 f^i_x\,\delta^{ij}
\delta(x-y)+A^{ij}\delta(x-y),\label{SZ-10a}\\
&&\omega^{ij}_2=g^i\,\delta^{ij}\delta'(x-y)+\frac12 g^i_x\,\delta^{ij}
\delta(x-y)+B^{ij}\delta(x-y).\label{SZ-10b}
\end{eqnarray}
Here $f^i=f^i(u^1,\dots,u^n),\ g^i=u^i\,f^i,\ f^i_x=\p_x f^i,\ g^i_x=\p_x g^i$, and
\begin{equation}\label{SZ-10c}
A^{ij}=\frac12\left(\frac{f^i}{f^j}f^j_iu^{j}_{x}-
\frac{f^j}{f^i}f^i_ju^{i}_{x}\right),
B^{ij}=\frac12\left(\frac{u^i f^i}{f^j}f^j_iu^{j}_{x}-
\frac{u^j f^j}{f^i}f^i_ju^{i}_{x}\right)
\end{equation}
where $f^a_b=\frac{\p f^a}{\p u^b}$.

Denote by ${\hat \Omega}^0$ the space of local functionals of the form
$$
{\bar f}=\int f(u,u_x,\dots,u^{(N)}) dx
$$
where $f$ is a smooth function of all of its arguments. Define
\begin{eqnarray}
&&{\hat H}^2(\LM;\omega_1,\omega_2)=\oplus_{m\ge 1} {\hat H}^2_m,\nn\\
&&{\hat H}^2_m=H^2_m(\LM;\omega_1,\omega_2)\cap
(d_1\hat{\Omega}^0\oplus d_2\hat{\Omega}^0).
\end{eqnarray}
Then ${\hat H}^2$ is the space of equivalence classes of infinitesimal quasitrivial deformations
of the bihamiltonian structure $(\omega_1, \omega_2)$.

\begin{thm}\label{mthm}
We have ${\hat H}^2_m=0$ for $m=1,3,4,\dots$ and
\begin{equation}\label{fpc}
{\hat H}_2^2=\{\sum_{i=1}^n \left(d_2\int (c_i(u^i)u^i_x\log u^i_x)dx
-d_1\int (u^ic_i(u^i)u^i_x\log u^i_x)dx\right)\}.
\end{equation}
Here $d_1, d_2$ are the differentials defined by the Poisson structures $\omega_1$ and $\omega_2$
respectively, $c_i(u^i)$ are arbitrary smooth functions of $u^i$. Moreover, two sets of
functions $\{c_i\}$ and $\{\tilde{c}_i\}$ define the same element in
${\hat H}^2$ iff $c_i=\tilde{c}_i$.
\end{thm}
We will use the symbol
$$A(u,u_x,\dots,u^{(N)}) \sim B(u,u_x,\dots,u^{(N)})$$
to indicate that the difference of the functions $A$ and $B$ is a
differential polynomial. In order to prove the above theorem we first need to prove some lemmas.
\begin{lem}\label{lemma-1}
Let $X=d_2 I-d_1 J\in \hat{H}^2$ with
$$I=\int G(u,u_x,\dots,u^{(N)}) dx,
\ J=\int H(u,u_x,\dots,u^{(N)}) dx,\ N\ge2.$$
Then the densities $G, H$ can be chosen to have the form
\begin{eqnarray}
G&\sim&\sum_{i=1}^n \frac{(u^{i,N})^2}{u^i_x}
P^i(u;u_x,\dots,u^{(N-2)};u^{i,N-1})+Q(u,\dots,u^{(N-1)}),\label{z6}\\
H&\sim&\sum_{i=1}^n \frac{(u^{i,N})^2}{u^i_x} u^i
P^i(u;u_x,\dots,u^{(N-2)};u^{i,N-1})+R(u,\dots,u^{(N-1)}).\label{z7}
\end{eqnarray}
Here $P^i(u;u_x,\dots,u^{(N-2)};u^{i,N-1})$ are differential polynomials, $Q, R$ are smooth functions,
and any nonzero differential polynomial $P^i(u;u_x,\dots,u^{(N-2)};u^{i,N-1})$ is indivisible by $u^i_x$,
\end{lem}

\begin{prf}
Denote by $X^i,i=1,\dots,n$ the components of the local vector field  $X$, from our
assumption we know that they are differential polynomials. We are to use this property
repeatedly to prove the lemma. Let us start with the polynomiality of
$\frac{\p X^i}{\p u^{j,2 N+1}}$. Denote
$$X^i_{j,m}=\frac{\p X^i}{\p u^{j,m}},\
G_{i,p;j,q}=\frac{\p^2 G}{\p u^{i,p}\p u^{j,q}},\
H_{i,p;j,q}=\frac{\p^2 H}{\p u^{i,p}\p u^{j,q}}.$$
By using the simple identity
$$\frac{\p}{\p u^{i,k}} \p_x^m=\sum_{l=0}^m \binom{m}{l}
\p_x^l \frac{\p}{\p u^{i,k-m+l}}$$
and the form (\ref{SZ-10a}), (\ref{SZ-10b}) of the bihamiltonian structure
$(\omega_1, \omega_2)$ we obtain the following formulae:
\begin{equation}\label{SZ-12}
(-1)^N X^i_{j,2 N+1}=g^i G_{i,N;j,N}-f^i H_{i,N;j,N}.
\end{equation}
It follows that the functions $G$ and $H$ satisfy the relations
\begin{equation}\label{SZ-15}
u^i G_{i,N;j,N}-H_{i,N;j,N}\sim0,\ (u^i-u^j)G_{i,N;j,N}\sim 0.
\end{equation}
So there exist smooth functions $a_i, b_i, c$,
such that
\begin{eqnarray*}
G&\sim&\sum_{i=1}^n a_i(u,\dots,u^{(N-1)},u^{i,N}),\\
H&\sim&\sum_{i=1}^n (u^i a_i(u,\dots,u^{(N-1)},u^{i,N})+
b_i(u,\dots,u^{(N-1)}) u^{i,N})+c(u,\dots,u^{(N-1)}).
\end{eqnarray*}
By substituting these expressions into the relations
$(-1)^N \frac{\p X^i}{\p u^{i,2N}}\sim 0$
we obtain
$$-(N+\frac12) f^i\,u^i_x\, \frac{\p^2 a_i}{\p u^{i,N}\p u^{i,N}}\sim0.$$
Thus we can find differential polynomials $p_i(u,\dots,u^{(N-1)},u^{i,N})$ and smooth functions $q_i(u,\dots,u^{(N-1)}),\,r_i(u,\dots,u^{(N-1)})$ such that
$$a_i=\frac{p_i(u,\dots,u^{(N-1)},u^{i,N})}{u^i_x}+q_i(u,\dots,u^{(N-1)}) u^{i,N}
+r_i(u,\dots,u^{(N-1)}).$$
Now the functions $G, H$ can be written in the form
\begin{eqnarray}
G&\sim&\sum_{i=1}^n \left(\frac{p_i(u,\dots,u^{(N-1)},u^{i,N})}{u^i_x}
+q_i(u,\dots,u^{(N-1)}) u^{i,N}\right)\nn\\
&&\quad+r(u,\dots,u^{(N-1)}),\label{SZ-16}\\
H&\sim&\sum_{i=1}^n \left(u^i\frac{p_i(u,\dots,u^{(N-1)},u^{i,N})}{u^i_x}
+s_i(u,\dots,u^{(N-1)}) u^{i,N}\right)\nn\\
&&\quad+e(u,\dots,u^{(N-1)}).\label{SZ-17}
\end{eqnarray}
Here $s_i, e$ are some smooth functions.
{\em In the above expression of $G, H$, we assume that the differential
polynomials $p_i$ do not contain terms that are linear and constant with
respect to $u^{i,N}$, such terms can be absorbed into the functions
$q_i\, u^{i,N},\, s_i\, u^{i,N}$ and  $r, e$}.

Assuming the form (\ref{SZ-16}) and (\ref{SZ-17}) of the functions $G, H$
we continue to use the polynomoality of $(-1)^N \frac{\p X^i}{\p u^{j, 2N}}$
with $i\ne j$ to obtain
$$u^i\left(G_{i,N;j,N-1}-G_{j,N;i,N-1}\right)-
\left(H_{i,N;j,N-1}-H_{j,N;i,N-1}\right)\sim0.$$
From these relations it follows that for indices $i\ne j$ we have
\begin{eqnarray}\label{SZ-18}
&&H_{i,N;j,N-1}-H_{j,N;i,N-1}\sim0,\ G_{i,N;j,N-1}-G_{j,N;i,N-1}\sim0\\
&&G_{i,N;i,N;j,N-1}-G_{i,N;i,N-1;j,N}\sim G_{i,N;i,N;j,N-1}\sim0\label{SZ-18-b}
\end{eqnarray}
The relation (\ref{SZ-18-b}) shows that we can adjust the differential
polynomials $p_i$ so that they have the form
$$p_i=p_i(u,\dots,u^{(N-2)},u^{i,N-1},u^{i,N}),\quad i=1,\dots,n.$$
Now by substituting the expression (\ref{SZ-17}) for the function $H$ into the
first relation of (\ref{SZ-18}) we arrive at
$$\frac{\p s_i}{\p u^{j,N-1}}-\frac{\p s_j}{\p u^{i,N-1}}\sim0,$$
by using the Poincar\'e lemma we can find differential polynomials
${\hat s}_1,\dots {\hat s}_n$ such that the identity
$$\frac{\p (s_i-{\hat s}_i)}{\p u^{j,N-1}}-\frac{\p (s_j-{\hat s}_j)}{\p u^{i,N-1}}=0$$
hold true. This identity implies the existence of a function $W(u,\dots,u^{(N-1)})$ satisfying
$$s_i \sim \frac{\p W}{\p u^{i,N-1}},\quad i=1,\dots,n.$$
So by adjusting the density $H$ of the functional $J$ to $H-\p_x W$ we can
assume that in the expression (\ref{SZ-17}) for the function $H$ the second term
$\sum_{i=1}^n s_i u^{i,N}$ does not appear. In a similar way, we can also assume that
the term $\sum_{i=1}^n q_i u^{i,N}$ in the expression (\ref{SZ-16})
of the density of the functional $I$ vanishes.

Finally, the relation $(-1)^N \frac{\p^2 X^i}{\p u^{i,2N-1}\p u^{i,N}}\sim 0$
implies that
\begin{equation}
\frac{N^2} 2 \frac{f^i\,u^i_{xx}}{u^i_x}\,
\frac{\p^3 p_i}{\p u^{i,N}\p u^{i,N}\p u^{i,N}}\sim 0, \quad i=1,\dots,n.
\end{equation}
So we can adjust the densities $G, H$
of the functionals $I, J$  so that they have the forms (\ref{z6}), (\ref{z7}).
The lemma is proved.
\end{prf}

Let us introduce the operators
$$Z^m_{ij}=\sum_{p\ge m}(-1)^p\binom{p}{m}
\frac{\p^2}{\p u^{i,p}\,\p u^{j,2N+m-p}},\ 1 \le i,j \le n,\ m \ge0.$$
It is easy to verify that these operators satisfy the identities $[\p_x,Z^{m}_{ij}]=Z^{m-1}_{ij}$ and,
moreover, we have the following lemma:

\begin{lem}\label{lemma-2}
For a functional $I=\int G(u,u^{(1)},\dots)dx$, denote
\begin{equation}\label{z7-b}
I_k=\frac{\delta I}{\delta u^k},\ k=1,\dots,n.
\end{equation}
Then for any indices $i,j,m$,
the following formulae hold true
$$Z^m_{ij}I_k=\sum_{s\ge0}\binom{s+m}{s}(-\p_x)^s
\frac{\p}{\p u^{k,s+m}}\left(\frac{\p I_i}{\p u^{j,2N}}\right)$$
\end{lem}
\begin{prf}
It is well known from the theory of  variational calculus that for any functional
$I$ we have the following identities:
$$\frac{\p}{\p u^{i,p}}\left(\frac{\delta I}{\delta u^k}\right)
=\sum_{t\ge p}(-1)^t \binom{t}{p} \p_x^{t-p}
\frac{\p}{\p u^{k,t}}\left(\frac{\delta I}{\delta u^i}\right)$$
From which it follows that
$$\frac{\p^2 I_k}{\p u^{i,p}\p u^{j,2N+m-p}}=\sum_{s \ge0}
\sum_{t\ge p} (-1)^{s+t} \binom{s+t}{p}\binom{s+t-p}{s} \p_x^s
\frac{\p^2 I_i}{\p u^{k,s+t}\p u^{j,2N+m-t}}.$$
By using this identity we obtain
\begin{eqnarray*}
&&Z^m_{ij} I_k=\sum_{p\ge0}(-1)^p\binom{p}{m}
\frac{\p^2 I_k}{\p u^{i,p}\,\p u^{j,2N+m-p}}\\
&=&\sum_{p\ge0}(-1)^p \binom{p}{m}
\sum_{s\ge0}\sum_{t\ge p} (-1)^{s+t} \binom{s+t}{p}\binom{s+t-p}{s}
\p_x^s \frac{\p^2 I_i}{\p u^{k,s+t}\p u^{j,2N+m-t}}\\
&=&\sum_{s\ge0}(-\p_x)^s \sum_{t\ge0}(-1)^t\binom{s+t}{s}
\left[\sum_{p=0}^t (-1)^{p} \binom{p}{m}\binom{t}{p}\right]
\frac{\p^2 I_i}{\p u^{k,s+t}\p u^{j,2N+m-t}}\\
&=&\sum_{s\ge0}(-\p_x)^s\binom{s+m}{s}\frac{\p^2 I_i}{\p u^{k,s+m}\p u^{j,2N}}.
\end{eqnarray*}
Here we assumed $\binom{p}{m}=0$ when $p\le m-1$ and we used the identity
$$\sum_{p=0}^t (-1)^{p} \binom{p}{m}\binom{t}{p}=(-1)^t \delta_{tm}.$$
The lemma is proved.
\end{prf}

\begin{lem}\label{lemma-3}
The polynomials $P^i$ defined in Lemma \ref{lemma-1} must vanish.
\end{lem}
\begin{prf}
Let $m$ be the highest order of the $x$-derivatives of $u^1,\dots u^n$ that
appear in the polynomials $P^i$. We first prove, by using
the polynomiality of $Z^{m-1}_{ij}X^k$, that $m$ must less than $3$.
To this end, let's assume at the moment that $m\ge 3$.
From the form (\ref{SZ-10a}), (\ref{SZ-10b}) of the bihamiltonian structure $(\omega_1, \omega_2)$
we know that the components of the vector field $X=d_2 I-d_1 J$ can be expressed as
\begin{eqnarray*}
X^k&=&
g^k\p_x\frac{\delta I}{\delta u^k}+
\frac{\p_x g^k}2\frac{\delta I}{\delta u^k}+
\sum_{\alpha=1}^n B^{k\alpha}\frac{\delta I}{\delta u^{\alpha}}\\
&&\qquad-f^k\p_x\frac{\delta J}{\delta u^k}-
\frac{\p_x f^k}2\frac{\delta J}{\delta u^k}-
\sum_{\alpha=1}^n A^{k\alpha}\frac{\delta J}{\delta u^{\alpha}}\\
\end{eqnarray*}
Since the highest order of the $x$-derivatives of $u^p$ that appear in
$\frac{\delta I}{\delta u^k}$ is $2N$, we have
\begin{eqnarray*}
&&Z^{m-1}_{ij}X^k \\
&=&g^k(\p_x\,Z^{m-1}_{ij}-Z^{m-2}_{ij})I_k+\frac{\p_x g^k}2Z^{m-1}_{ij}I_k
+\sum_{\alpha=1}^nB^{k\alpha}Z^{m-1}_{ij}I_{\alpha}\\
&& \qquad -f^k(\p_x\,Z^{m-1}_{ij}-Z^{m-2}_{ij})J_k-\frac{\p_x f^k}2
Z^{m-1}_{ij}J_k-\sum_{\alpha=1}^nA^{k\alpha}Z^{m-1}_{ij}J_{\alpha}
\end{eqnarray*}
Here $I_k, J_k$ are defined as in (\ref{z7}).
By using Lemma \ref{lemma-1} and \ref{lemma-2} we know that
\begin{equation*}
\frac{\p I_i}{\p u^{j,2N}}\sim (-1)^N\frac{2 P^i}{u^i_x}\delta_{ij}.
\end{equation*}
and
\begin{eqnarray*}
&&Z^{m-1}_{ij}I_k\sim \left(\frac{\p}{\p u^{k,m-1}}-m\p_x\frac{\p}{\p u^{k,m}}\right)
\frac{\p I_i}{\p u^{j,2N}},
\\
&&Z^{m-2}_{ij}I_k\sim \left(\frac{\p}{\p u^{k,m-2}}-(m-1)\p_x\frac{\p}{\p u^{k,m-1}}
+\frac{m(m-1)}2\p_x^2\frac{\p}{\p u^{k,m}}\right)
\frac{\p I_i}{\p u^{j,2N}}.
\end{eqnarray*}
We can get similar expression for $Z^{m-1}_{ij}J_k$ and $Z^{m-2}_{ij}J_k$.
By using these formulae, we see that for the case $i=j\ne k$
the term with the highest
power of $\frac1{u^i_x}$ in the expression of $Z^{m-1}_{ij}X^k$ is given by
\begin{equation}\label{ux3}
(-1)^{N}\ 2m(m+1)f^k(u^i-u^k)\frac{(u^i_{xx})^2}{(u^i_x)^3}
\frac{\p P^i}{\p u^{k,m}}
\end{equation}
From the fact that $P^i$ is indivisible by $u^i_x$ and $Z^{m-1}_{ii}X^k$ is a differential polynomial
it follows that $P^i$ does not depend on $u^{k,m}$ for $k \ne i$. In the case when  $i=j=k$ we have
\begin{equation}\label{ux1}
Z^{m-1}_{ij}X^k \sim (-1)^{N+1} m^2 f^i\frac{u^i_{xx}}{u^i_x}
\frac{\p P^i}{\p u^{i,m}}.
\end{equation}
So $P^i$ does not depend on $u^{i,m}$ either. Thus we proved that the highest order $m$ of
the $x$-derivatives of $u^1,\dots u^n$
that appear in the polynomial $P^i$ must less than $3$. To complete the proof of the
lemma we use the polynomiality of $Z^1_{ij}X^k$. In the same way as we did above,
we can prove that the terms (\ref{ux3})
for the case of $m=2$ is a differential polynomial, so $P^i$
does not depend on $u^k_{xx}$ for $i \ne k$. Then the counterpart of (\ref{ux1}) for the case of
$m=2$ has the form
\begin{equation}
Z^1_{ii}X^i \sim \frac{(-1)^{N+1} f^i}{u^i_x}
\left(4u^i_{xx}\frac{\p P^i}{\p u^i_{xx}}+(2N-2)P^i\right)
\end{equation}
which implies $P^i=0$. The lemma is proved.
\end{prf}

Now we can prove the main result of this section.

\begin{prfM}
By using the above lemma, we know that for any element of $\hat{H}^2$ we can
choose its representative
$X\in \Ker(d_1 d_2)$ of the form
\begin{equation}
X=d_2 I-d_1 J,\quad I=\int G(u,u_x) dx, \ J=\int H(u,u_x) dx.
\end{equation}
Then the polynomiality of
\begin{equation}
\frac{\p X^i}{\p u^{j,3}}=f^i \frac{\p^2 H}{\p u^i_x \p u^j_x}-g^i \frac{\p^2 G}{\p u^i_x \p u^j_x}
\end{equation}
allows us to adjust the vector field $X$ such that the functions $G$ and
$H$ have the expression
\begin{equation}
G=\sum_{i=1}^n h_i(u^1,\dots,u^n, u^i_x),\quad
H=\sum_{i=1}^n u^i h_i(u^1,\dots,u^n, u^i_x).
\end{equation}
By using the identity
\begin{equation}
\frac{\p X^i}{\p u^i_{xx}}
=\frac32 f^i u^i_x \frac{\p^2 h_i}{\p u^i_x \p u^i_x}
\end{equation}
we see that the functions $h_i$ must take the form
\begin{equation}
h_i=c_i(u) u^i_x \log u^i_x+{\rm {differential\, polynomial}}.
\end{equation}
Now from the explicit form of  $\frac{\p X^i}{\p u^j_{xx}}$
we know that
\begin{equation}
(u^i-u^j) \frac{\p c_j}{\p u^i} \log u^j_x
\end{equation}
are differential polynomials, thus we have $\frac{\p c_j}{\p u^i}=0$ for $i\ne j$,
and $c_i$ depend only on $u^i$. So we proved that any element of $\hat{H}^2$ has
a representative of the form given in the right hand side of (\ref{fpc}).

On the other hand, given any vector field $X$ with the form given in the right
hand side of (\ref{fpc}), we can easily verify that its components have the
expressions
\begin{equation}\label{SZ-22}
X^i=\sum_{j=1}^n\left[ \left(\frac12 \delta_{ij} \p_x f^i+
A^{ij}\right) c_j u^j_x+(2 \delta_{ij} f^i
-L^{ij}) \p_x \left(c_j u^j_x\right)\right].
\end{equation}
Here
\begin{equation}
L^{ij}=\frac12 \delta_{ij} f^i+\frac{(u^i-u^j) f^i}{ 2 f^j}\,
\frac{\p f^j}{\p u^i}.\end{equation}
It shows that $X^i$ are differential polynomials and thus $X$ is a
representative of an element of $\hat{H}^2$.

Finally, we are left to show that a
vector field $X$ of the form given in the right hand side of (\ref{fpc}) is
trivial if and only if $c_1=\dots =c_n=0$. From the expression (\ref{SZ-22})
it follows that the triviality of the  vector field $X$ is equivalent to the
existence of functions $\alpha_i(u), \beta_i(u),\, i=1,\dots, n$ such that
the vector fields $X$ can be expressed as  ${\tilde X}=d_2{\tilde I}-d_1 {\tilde J}$,
where the functionals ${\tilde I}$ and ${\tilde J}$ have the form
\begin{equation}
{\tilde I}=\int \sum_{i=1}^n \alpha_i(u) u^i_x dx,
\quad {\tilde J}=\int \sum_{i=1}^n \beta_i(u) u^i_x dx.
\end{equation}
The coefficient of $u^i_{xx}$ of the $i$-th component of $X$ is given by
$2 f^i c_i$, while that of ${\tilde X}$ equals zero.
Thus we must have $c_i=0,\, i=1,\dots, n$. The theorem is proved.
\end{prfM}

\noindent {\it Proof of Theorem \ref{mthm-z1} and Theorem \ref{mthm-z2}} \
Let us assume that the
hydrodynamic bihamiltonian structure  $(\omega_1,\omega_2)$
has two $N$-th order quasitrivial
deformations of the form
\begin{eqnarray}
&&(\omega_1,\omega_2+\sum_{m=1}^N \ve^m P_m)+{\mathcal O}(\ve^{N+1}),\label{z15-a}\\
&&(\omega_1,\omega_2+\sum_{m=1}^N \ve^m P_m+\ve^N Q)+{\mathcal O}(\ve^{N+1}).\label{z15-b}
\end{eqnarray}
Here $P_m\in \Omega^2_{m+2}, Q\in \Omega_{N+2}^2$. Due to our assumption, we can find a quasi-Miura
transformation of the form (\ref{qmiura}) that transforms the bihamiltonian structure (\ref{z15-a})
to $(\omega_1,\omega_2)+{\mathcal O}(\ve^{N+1})$. Then this same quasi-Miura transformation transforms the
bihamiltonian structure (\ref{z15-b}) to
\begin{equation}
(\omega_1,\omega_2+\ve^N Q)+{\mathcal O}(\ve^{N+1}).\label{z16}
\end{equation}
It is also a quasitrivial deformation of the bihamiltonian structure $(\omega_1,\omega_2)$, so we are able to
find a quasi-Miura transformation
that transforms $(\omega_1,\omega_2)$ to (\ref{z16}). Such a quasi-Miura transformation can be represented
by some vector fields $Y_1,\dots,Y_N$ in the form
\begin{eqnarray}
&&\omega_1=e^{-\ve^N {ad}_{{Y}_N}}\dots e^{-\ve {ad}_{{Y}_1}} \omega_1+{\mathcal O}(\ve^{N+1}),\nn\\
&&
\omega_2+\ve^N Q=e^{-\ve^N {ad}_{Y_N}}\dots e^{-\ve {ad}_{{Y}_1}}
\omega_2+{\mathcal O}(\ve^{s+1}).\label{z5-b}
\end{eqnarray}
From the above identities it follows that
$d_1 Y_N=0,\ Q+d_2 Y_N=0$, so there exists a functional $I$ such that
$Y_N=d_1 I,\ Q=d_1d_2 I$. On the other hand,
the compatibility of $(\omega_1,\omega_2+\ve^N Q)+{\mathcal O}(\ve^{N+1})$ implies the existence
of a vector field $X\in \Omega^1_{N}$ satisfying $Q=d_1 X$. From the above two expressions of $Q$
we see that we can express the vector field $X$ as
$$X=d_2I-d_1J$$
with certain functional $J\in{\hat \Omega}^{0}$.

Now the results of Theorem \ref{mthm} lead to the following conclusions:\newline
\noindent{\em 1.} If $N \ne 2$, then $I$ and $J$ must be diffrential polynomials, so the two deformations (\ref{z15-a})
and (\ref{z15-b}) are related by a Miura transformation
\begin{equation}
u\mapsto u-\ve^N d_1d_2 I.
\end{equation}
Theorem \ref{mthm-z1} is proved.\newline
\noindent{\em 2.} Any second order deformation $(\omega_1, \omega_2+\ve P_1+\ve^2 P_2)+{\mathcal O}(\ve^3)$
is equivalent to a second order deformation of the form $(\omega_1, \omega_2+\ve^2 {\tilde P}_2)
+{\mathcal O}(\ve^3)$. By applying the results of Theorem \ref{mthm} to the case with $N=2$,
we see that modulo a Miura transformation the deformed bihamiltonian structure can be represented
in the form
\begin{equation}
(\omega_1,\omega_2+\ve^2 d_1(d_2I-d_1J))+{\mathcal O}(\ve^3) \label{normal-a}
\end{equation}
for some functionals $I, J$ defined by
\begin{equation}
I=\int \sum_{i=1}^n c_i(u^i)u^i_x\log u^i_x dx,\quad
J=\int u^i c_i(u^i) u^i_x \log u^i_x dx.\label{normal-b}
\end{equation}
On the other hand, it is easy to see that any functionals $I, J$ of the above form define
a second order quasitrivial deformation of the bihamiltonian $(\omega_1, \omega_2)$.
Theorem \ref{mthm-z2} is proved.\hfill $\Box$

From the proof of the main theorems it follows that the any equivalence class of quasitrivial
deformations of the bihamiltonian structure $(\omega_1, \omega_2)$ has a unique representative
of the form (\ref{normal-a}), (\ref{normal-b}) which corresponds to an element of the modified
cohomology ${\hat H}^2$.

\section{Some examples}\label{sec-5}
In this section, we consider as examples the deformations of the
bihamiltonian structures of hydrodynamic type that are related to the KdV and
the nonlinear Schr\"odinger equations, these deformations yield the the bihamiltonian
structures for the Camassa-Holm hierarchy \cite{CH,CHH,Fokas,Fu, FF} and its generalization.

Let us first consider deformations of the bihamiltonian structure (\ref{z2}). The
class of deformations that corresponds to the element of ${\hat H}^2$ (see Theorem \ref{mthm})
with
$
c(u)=-\frac1{24}
$
has a representative
\begin{eqnarray}
&&\{u(x),u(y)\}_1=\delta'(x-y),\nn\\
&&\{u(x),u(y)\}_2=u(x)\delta'(x-y)+\frac12\,u(x)'\delta(x-y)+\frac{\ve^2}8 \delta'''(x-y).\label{z8}
\end{eqnarray}
Here we redenote $u^1=u, c_1(u)=c(u)$. It is just the well known bihamiltonian structure for the KdV
hierarchy \cite{gardner,magri,ZF}. Now if we take
$
c(u)=-\frac1{24} u,
$
then the corresponding class of deformations has the following representative
\begin{eqnarray}
&&\{u(x),u(y)\}_1=\delta'(x-y)-\frac{\ve^2}8 \delta'''(x-y),\nn\\
&&\{u(x),u(y)\}_2=u(x)\delta'(x-y)+\frac12\,u(x)'\delta(x-y).\label{z9}
\end{eqnarray}
In fact, it is equivalent to the bihamiltonian structure
\begin{equation}
(\omega_1,\omega_2+\ve^2 d_1 (d_2 I-d_1 J))+{\mathcal O}(\ve^3)
\end{equation}
under the Miura transformation
$$
u\mapsto u+\frac{\ve^2}{16} u''.
$$
Here $(\omega_1,\omega_2)$ denotes the bihamiltonian structure (\ref{z2}) and
the functionals $I$ and $J$ are defined by
$$
I=-\frac{1}{24} \int u(x) u'(x) \log{u'(x)} dx,\quad J=-\frac{1}{24} \int u(x)^2 u'(x) \log{u'(x)} dx.
$$

The related bihamiltonian hierarchy of integrable systems is the
Camassa-Holm hierarchy that is well known in soliton theory.
It can be expressed by the following bihamiltonian
recursion relations:
\begin{equation}\label{1-ch-h}
\frac{\pal u}{\pal t^q}=\{u(x), H_q\}_1=\frac2{2 q+1} \{u(x),H_{q-1}\}_2,\quad q\ge 0.
\end{equation}
Here we start from the Casimir $H_{-1}=\int u(x) dx$ of the first Poisson bracket, and then
determine the Hamiltonians $H_q, q\ge 0$ recursively from the above relation. The recursive
procedure of finding the Hamiltonians $H_q$ is guaranteed by the triviality of the first Poisson
cohomology of the Poisson structure $\omega_1$ \cite{magri2,DZ1, get}.
The first nontrivial flow $\frac{\pal}{\pal t}=\frac{\pal}{\pal t^1}$ of the hierarchy can be put into the form
\begin{equation}\label{1-ch}
(v-\frac{\ve^2}{8} v_{xx})_{t}=v v_x-\frac{\ve^2}{12} v_x v_{xx}-\frac{\ve^2}{24} v v_{xxx}.
\end{equation}
Here the dependent variable $v$ is defined by
\begin{equation}
u=v-\frac{\ve^2} 8 v_{xx}.
\end{equation}
If we change the time variable as $t^1\mapsto t=-\frac13 t^1$ and put $\ve^2=8$, then the
resulting equation is just the
Camassa-Holm shallow water wave equation \cite{CH,CHH,Fokas,Fu,FF},
which possesses most of the important properties
of an integrable system. In particular, it has the following Lax pair representation
\begin{eqnarray}
\ve^2 \phi_{xx}&=&\left(2-\frac{8 v-\ve^2 v_{xx}}{2\lambda}\right)\phi,\label{lax-ch-a}\\
\phi_{t}&=&\frac13 (\lambda+v)\phi_x-\frac{v_x}6\phi\label{lax-ch-b}
\end{eqnarray}
and its initial value problems can be solved by using the inverse scattering method. The Camassa-Holm
equation also possesses some features that are distinguished from the usual
KdV-type integrable systems, such as
the existence of peaked solitons, the nonlinear dependence of the arguments of its
algebraic-geometric solutions on the spatial variable $x$ \cite{Alb}
and the non-existence of tau function\cite{DZ1}. We will call the equation (\ref{1-ch}) and the hierarchy
(\ref{1-ch-h}) the Camassa-Holm equation and the Camassa-Holm hierarchy respectively.

The quasitriviality of the bihamiltonian structure (\ref{z8}), (\ref{z9}) can be deduced from a result of \cite{DZ1}
on the quasitriviality of a general class of bihamiltonian structures. Details on this aspect will
be given in a subsequent publication.

For the choice of a general smooth function $c(u)$, we do not have at this moment an explicit
expression of the correspondent class of deformations of the bihamiltonian structure (\ref{z2}).
At the approximation up to $\ve^4$ Lorenzoni obtained the expression of a representative of
the corresponding class of deformations, and we can in fact go further to show that his result can be
modified to reach the approximation up to higher orders of $\ve$. This fact strongly indicates the
existence of a full deformation of the bihamiltonian structure (\ref{z2}) for any smooth function $c(u)$,
or equivalently, to the vanishing of the third bihamiltonian cohomologies
$H^3_m({\mathcal L};\omega_1,\omega_2),\ m\ge 5$ of the
bihamiltonian structure (\ref{z2}).

We now consider the deformations of the following bihamiltonian structure
\begin{eqnarray}
&&\{w_1(x), w_1(y)\}_1=\{w_2(x), w_2(y)\}_1=0,\nn\\
&&\{w_1(x), w_2(y)\}_1=\delta'(x-y).\label{z10-a}\\
&&\{w_1(x), w_1(y)\}_2=2 \delta'(x-y),\nn\\
&&\{w_1(x), w_2(y))\}_2=w_1(x) \delta'(x-y)+w_1'(x) \delta(x-y),\nn\\
&&\{w_2(x), w_2(y)\}_2=\left[w_2(x)\pal_x+\pal_x w_2(x))\right] \delta(x-y).\label{z10}
\end{eqnarray}
It is related to the Frobenius manifold with potential \cite{CDZ}
$$
F=\frac12 w_1^2 w_2+\frac12 w_2^2\left(\log w_2-\frac32\right).
$$
The canonical coordinates of this bihamiltonian structure are given by
\begin{equation}
u^{1,2}=w_1\pm2{\sqrt{w_2}}\ .
\end{equation}
Let us consider the following two classes of deformations:\par
\vskip 0.1truecm
\noindent{\em Case 1}. We take the element of ${\hat H}^2$
with
$
c_1(u)=c_2(u)=-\frac1{24},
$
then the corresponding class of deformations has a representative
\begin{eqnarray}
&&\{w_1(x), w_1(y)\}_1=\{w_2(x), w_2(y)\}_1=0,\nn\\
&&\{w_1(x), w_2(y)\}_1=\delta'(x-y).\label{z14-a}\\
&&\{w_1(x), w_1(y)\}_2=2 \delta'(x-y),\nn\\
&&\{w_1(x), w_2(y))\}_2=w_1(x) \delta'(x-y)+w_1'(x) \delta(x-y)-\ve
\delta''(X-Y),\nn\\
&&\{w_2(x), w_2(y)\}_2=\left[w_2(x)\pal_x+\pal_x w_2(x))\right] \delta(x-y).\label{z14-b}
\end{eqnarray}
To see this, let us denote by $\omega_1, \omega_2$ the two bivectors of the bihamiltonian
structure (\ref{z10-a}), (\ref{z10}), and by $I, J$ the functionals
\begin{eqnarray}
&&I=-\int \frac1{24}\left( (u^1_x \log u^1_x+u^2_x \log u^2_x\right) dx,\nn\\
&&J=-\int \frac1{24}\left( (u^1 u^1_x \log u^1_x+u^2 u^2_x \log u^2_x\right) dx,
\end{eqnarray}
then by a direct computation it can be verified that
the bihamiltonian structure (\ref{z14-a}), (\ref{z14-b}) is equivalent to the bihamiltonian structure
\begin{equation}
(\omega_1,\omega_2+\ve^2 d_1 (d_2 I-d_1 J))+{\mathcal{O}}(\ve^3)
\end{equation}
under the Miura transformation
\begin{eqnarray}
&&w_1\mapsto w_1+\frac{\ve}{2\sqrt{3}}\,\frac{w_{2,x}}{w_2}+
\ve^2\left(\frac1{12}-\frac1{4\sqrt{3}}\right)\left(\frac{w_{1,xx}}{w_2}-\frac{w_{1,x} w_{2,x}}{w_2^2}\right),\nn\\
&&w_2\mapsto w_2+\ve\left(-\frac1{2}+\frac1{2\sqrt{3}}\right) w_{1,x}.
\end{eqnarray}
The bihmailtonian hierarchy of integrable systems that is related to this bihamiltonian
structure is called the {\em extended NLS hierarchy}, the algebraic properties of this
hierarchy together with its relation to the $CP^1$ topological sigma model is studied in
detail in \cite{CDZ,DZ3}. It is also shown in \cite{CDZ} that this hierarchy is equivalent to
the extended Toda hierarchy \cite{get2,YZ} which contains the standard Toda lattice hierarchy.\par
\vskip 0.1truecm
\noindent{\em Case 2}. Let us take the element of ${\hat H}^2$
with
$
c_1(u)=-\frac{(u^1)^2}{24}, c_2(u)=-\frac{(u^2)^2}{24},
$
then the correspondent class of deformations has a representative of the form
\begin{eqnarray}
&&\{w_1(x), w_1(y)\}_1=\{w_2(x), w_2(y)\}_1=0,\nn\\
&&\{w_1(x), w_2(y)\}_1=\delta'(x-y)-\ve
\delta''(x-y).\label{z11-a}\\
&&\{w_1(x), w_1(y)\}_2=2 \delta'(x-y),\nn\\
&&\{w_1(x), w_2(y))\}_2=w_1(x) \delta'(x-y)+w_1'(x) \delta(x-y),\nn\\
&&\{w_2(x), w_2(y)\}_2=\left[w_2(x)\pal_x+\pal_x w_2(x))\right] \delta(x-y).\label{z11-b}
\end{eqnarray}
Denote by $I, J$ the functionals
\begin{eqnarray}
&&I=-\int \frac1{24}\left( (u^1)^2 u^1_x \log u^1_x+(u^2)^2 u^2_x \log u^2_x\right) dx,\nn\\
&&J=-\int \frac1{24}\left( (u^1)^3 u^1_x \log u^1_x+(u^2)^3 u^2_x \log u^2_x\right) dx,
\end{eqnarray}
then it can be verified that
the bihamiltonian structure (\ref{z11-a}), (\ref{z11-b}) is equivalent to the bihamiltonian structure
$$
(\omega_1, \omega_2+\ve^2 d_1 (d_2 I-d_1 J))+{\mathcal O}(\ve^3)
$$
modulo a Miura transformation of the form
\begin{eqnarray*}
&&w_1\mapsto
w_1+\ve^2\left(\frac{w_1^2+4w_2}{24w_2}w_{1,x}\right)_x
+{\mathcal O}(\ve^3) \\
&&w_2\mapsto w_2+\ve \left(\frac{w_1^2}4-w_2\right)_x-
\ve^2\left(\left(\frac{w_1^2+4w_2}{24w_2}-1\right)w_{2,x}\right)_x
+{\mathcal O}(\ve^3)
\end{eqnarray*}
A hierarchy of integrable systems can be obtained by using the bihamiltonian recursion
relation
\begin{equation}
\{w_i(x), H_{q-1}\}_2=(q+1) \{w_i(x), H_q\}_1,\quad q\ge 0.
\end{equation}
Here we start from the Casimir $H_{-1}=\int w_2(x) dx$ of the first Poisson bracket, and then
determine the Hamiltonians $H_q, q\ge 0$ recursively by using the above relation.
The flows of the bihamiltonian hierarchy is then given by
\begin{equation}\label{2-ch-h}
\frac{\pal w_i}{\pal t^q}=\{w_i(x), H_q\}_1,\quad q\ge 0.
\end{equation}
The first flow $\frac{\pal}{\pal t_0}$ corresponds to the translation along the spatial
variable $x$, and the second flow $\frac{\pal}{\pal t}=\frac{\pal}{\pal t^1}$ has the form
\begin{eqnarray}
&&(\varphi_1-\ve \varphi_{1,x})_{t}=(\varphi_2+\frac12 \varphi_1^2-\frac{\ve}2 \varphi_1 \varphi_{1,x})_x, \\
&&(\varphi_2+\ve \varphi_{2,x})_{t}=(\varphi_1 \varphi_2+\frac{\ve} 2 \varphi_1 \varphi_{2,x})_x.\label{2-ch}
\end{eqnarray}
Here $\varphi_1, \varphi_2$ are defined by $ w_1=\varphi_1-\ve \varphi_{1,x},\
w_2=\varphi_2+\ve \varphi_{2,x}$. By introducing the new variables
$$
v_1=\varphi_1,\quad v_2=\varphi_2+\ve \varphi_{2,x}-\frac14 (\varphi_1-\ve \varphi_{1,x})^2
$$
we can rewrite the above system of equations in the following form
\begin{eqnarray}
&&(v_1-\ve^2 v_{1,xx})_t=\left(v_2+\frac34 v_1^2-\ve^2 (\frac12 v_1 v_{1,xx}+\frac14 v_{1,x}^2)\right)_x,\label{z12-a}\\
&&v_{2,t}=\frac12 v_1 v_{2,x}+v_2 v_{1,x}.\label{z12-b}
\end{eqnarray}
It easily follows from the above expression that the system of equations (\ref{z12-a}), (\ref{z12-b})
is reduced to the Camassa-Holm equation (\ref{1-ch}) under the constraint
\begin{equation}
v_2=0
\end{equation}
together with the rescaling $t\mapsto \frac32 t,\ \ve^2\mapsto \frac18{\ve^2}$. So we
can view the hierarchy (\ref{2-ch-h}) as a natural 2-component generalization of the Camassa-Holm hierarchy
(\ref{1-ch-h}). The following Lax pair formalism of the system (\ref{z12-a}), (\ref{z12-b})
manifests the above observation:
\begin{eqnarray}
\ve^2 \phi_{xx}&=&\left(\frac14-\frac{v_1-\ve^2 v_{1,xx}}{2\lambda}-
\frac{v_2}{\lambda^2}\right)\phi,\\
\phi_t&=&\frac12(\lambda+v_1)\phi_x-\frac{v_{1,x}}4\phi.
\end{eqnarray}
When we put $v_2=0$ this Lax pair is reduced to the one that is given in (\ref{lax-ch-a}), (\ref{lax-ch-b}).

The quasitriviality of the bihamiltonian structure
(\ref{z14-a}), (\ref{z14-b}) can be verified by using the method given in \cite{DZ1}. However, at this
moment we do not have a proof for the quasitriviality of the bihamiltonian structure (\ref{z11-a}), (\ref{z11-b}).
In order to use the approach of \cite{DZ1} to prove its quasitriviality we need to construct a bihamiltonian
hierarchy  of the form (\ref{2-ch-h}) that corresponds to the Casimir $\int w_1(x) dx$
of the first Poisson bracket, since this functional is also a Casimir of the second Poisson bracket,
the usual bihamiltonian recursion procedure fails to yield the needed Hamiltonians in a direct way.
We will consider in detail the propertities of the above 2-component Camassa-Holm hierarchy and its
further generalizations in a separate publication.

\section{Concluding remarks}
For any semisimple bihamiltonian structure of hydrodynamic type,
we classify its infinitesimal quasitrivial deformations. We show
that the equivalence classes of its second order quasitrivial
deformations are parameterized by $n$ arbitrary functions of one
variable, and we prove that any class of its quasitrivial
deformations is uniquely determined by its corresponding class of
second order deformations. We end this paper with the following
two remarks:\par \vskip 0.1truecm \noindent{\em Remark 1.} At a
first glance the condition of quasitriviality seems to be highly
non-trivial, however, a careful study shows that any deformation
of the semisimple bihamiltonian structure of the form
(\ref{SZ-10a}), (\ref{SZ-10b}) is quasitrivial at least for the
case of $n=1$, this fact together with the quasitriviality of any
tau-symmetric bihamiltonian structure \cite{DZ1} indicates the
validity of quasitriviality for any deformation of the semisimple
bihamiltonian structure of the form (\ref{SZ-10a}),
(\ref{SZ-10b}). An even more optimistic conjecture is the
existence of a full deformation of a semisimple bihamiltonian
structure of hydrodynamic type with a given second order
deformation. In the language of bihamiltonian
cohomology we can formulate the above conjectures as follows:
\begin{cjt}
For any semisimple bihamiltonian structure of hydrodynamic type $({\mathcal L}(M); \omega_1, \omega_2)$
we have ${H}^2({\mathcal L}(M);\omega_1,\omega_2)
={\hat H}^2({\mathcal L}(M);\omega_1,\omega_2)$,
and
the third bihamiltonian cohomologies $H^3_m({\mathcal L}(M); \omega_1,\omega_2)$ for $m\ge 5$ are trivial.
\end{cjt}
\noindent{\em Remark 2.} On the formal loop space of any
semisimple Frobenius manifold there is defined a semisimple
bihamiltonian structure of hydrodynamic type \cite{Du1}, a class
of deformations of such bihamiltonian structure was constructed in
\cite{DZ1}, these deformations correspond to the element of the
second cohomology ${\hat H}^2$ with $c_1=\dots=c_n=-\frac1{24}$,
they are compatible with the universal identities satisfied by the
Gromov-Witten invariants of smooth projective varieties, for this
reason we call them the topological deformations. The
corresponding bihamiltonian hierarchy of integrable systems
satisfies, in the sense of \cite{DZ1}, the properties of
tau-symmetry and linearization of the Virasoro symmetries. If we
drop the requirement of linearization of the Virasoro symmetries,
then the resulting tau symmetric bihamiltonian structure must
correspond to an element of the second cohomology ${\hat H}^2$
with constant $c_1(u)=c_1,\dots, c_n(u)=c_n$. An example of such
bihamiltonian structures is given by the one that is obtained by
using the Drinfeld-Sokolov construction for the affine Lie algebra
of type $B_2$ \cite{CP, DZ2,EYY}, in this case the corresponding
element of the second cohomology ${\hat H}^2$ is determined by the
constant functions $c_1=-\frac16, c_2=-\frac1{12}$.

\vskip 0.2truecm \noindent{\bf Acknowledgments.} The authors are
grateful to Boris Dubrovin for helpful suggestions and comments on
the research of the subject. The researches of Y.Z. were partially
supported by the Chinese National Science Fund for Distinguished
Young Scholars grant No.10025101 and the Special Funds of Chinese
Major Basic Research Project ``Nonlinear Sciences''.

\end{document}